\documentclass[aps,onecolumn,12pt,a4paper]{revtex4}
\usepackage{amsmath}
\usepackage{amssymb}
\newcommand{\beq}{\begin{equation}}
\newcommand{\eeq}{\end{equation}}
\newcommand{\bseq}{\begin{eqnarray}}
\newcommand{\eseq}{\end{eqnarray}}

\begin{document}
\title{\Large On Wilson's Theorem and Polignac Conjecture}
\author{Cong Lin and Li Zhipeng}
\affiliation{Hwa Chong Junior College, 661 Bukit Timah Road,
Singapore 269734}
\begin{abstract}
ABSTRACT. We introduce Wilson's Theorem and Clement's result and
present a necessary and sufficient condition for $p$ and $p+2k$ to
be primes where $k\in \mathbb{\mathbb{Z}^{+}}$. By using
Simionov's Theorem, a generalized form of Wilson's Theorem, we
derive improved version of Clement's result and characterization
of Polignac twin primes which parallels previous
characterizations. A straightforward method for determining the
coefficients in the derivation is also discussed.

\end{abstract}
\maketitle

\section{Introduction}
Prime numbers, whole numbers that are not the products of two
smaller numbers, are the bases of the arithmetics. In the
nineteenth century it was shown that the number of primes less
than or equal to $n$ approaches $n/(\log_e n)$ (as $n$ gets very
large); so a rough estimate for the n'th prime is $n \log_e n$.
The Sieve of Eratosthenes is still the most efficient way of
finding all very small primes (e.g., those less than 1,000,000).
However, most of the largest primes are found using special cases
of Lagrange's Theorem from group theory. In 1984 Samuel Yates
defined a titanic prime to be any prime with at least 1,000
digits.  When he introduced this term there were only 110 such
primes known; now there are over 1000 times that many! And as
computers and cryptology continually give new emphasis to search
forever-larger primes, this number will continue to grow.   The
problem of distinguishing prime numbers from composite numbers and
of resolving the latter into their prime factors is known to be
one of the most important and useful in arithmetic. There are
still many open questions  relating to prime numbers, among which
are the famous Twin Primes Conjecture and Goldbach's
Conjecture\cite{Goldbach}.

In relation to the unsolved conjectures, the properties of prime
numbers were studied and several theorems and corollaries
developed. In 1849, Alphonse de Polignac (1817-1890) made the
general conjecture that for any positive integer $k$, there are
infinitely many pairs of primes differ by $2k$. The case for which
$k=1$ is the twin primes case. Since then, numerous attempts to
prove this conjecture had been made. Based on heuristic
considerations, a law (the twin prime conjecture) was developed,
in 1922, by Godfrey Harold Hardy (1877-1947) and John Edensor
Littlewood (1885-1977) to estimate the density of twin primes\cite
{Hardy}. Twin prime characterization was discussed later by
Clement in 1949 \cite{Clement}. But an effective approach to the
conjecture remains undeveloped. It is thus worthwhile to
re-examine the derivation of Clement's result from a more
fundamental theorem in number theory, the Wilson's Theorem. By
generalizing and improving Wilson's Theorem and Clement's result,
a superior approach to Polignac Conjecture can be adapted to
obtain useful partial results to the mysterious problems.
\section{Wilson's Theorem and Clement's result}
Wilson derived his theorem on sufficient and necessary condition
for a number $p$ to be prime.  \bseq \textrm{\textbf{Theorem 1}}
\qquad \qquad \qquad (p-1)! \equiv -1 \pmod{p} \quad \textrm {iff
$p$ is prime.} & & \qquad \qquad \qquad \qquad \qquad \nonumber
\eseq In 1949, Clement \cite{Clement, frwebsite} formulated
another theorem based on Wilson's Theorem. The following is our
derivation of Clement's result,\bseq \textrm{\textbf{Theorem
2}}\qquad 4[(p-1)!+1] \equiv -p \pmod{p(p+2)} & & \textrm{ iff
$p$, $p+2$ are primes.} \qquad \qquad \qquad \nonumber \eseq
\textbf{Proof.}\quad \qquad Obviously, $p+2$ is not prime when
$p=2$. So we exclude $p=2$ from our discussion, i.e. $p$ and $p+2$
are odd prime numbers. By Wilson's Theorem, when $p \neq 2$, \bseq
& & (p-1)!
\equiv -1 \pmod{p} \qquad \textrm{ iff $p$ is odd prime.} \nonumber \\
& & \Rightarrow (p-1)! +1  \equiv  0 \pmod{p} \qquad \textrm{ iff
$p$ is odd prime.} \nonumber \\ & & \Rightarrow 4[(p-1)! +1]
\equiv 0
\pmod{p} \textrm{ and } p\neq4 \qquad \textrm{ iff $p$ is odd prime.} \nonumber \\
& & \Rightarrow 4[(p-1)! +1] \equiv  -p \pmod{p} \textrm{ and } p
\neq 4 \qquad \textrm{ iff $p$ is odd prime.}\eseq Similarly,
\bseq & & (p+1)! + 1 \equiv 0 \pmod{p+2} \textrm{ iff $p+2$ is odd
prime.} \nonumber
\\ & &
\Rightarrow (p^{2} + p)(p-1)! + 1  \equiv  0 \pmod{p+2} \qquad
\textrm{ iff $p+2$ is odd prime.} \nonumber
\\ & &
\Rightarrow [(p+2)(p-1)+2](p-1)! + 1  \equiv  0 \pmod{p+2} \qquad
\textrm{ iff $p+2$ is odd prime.} \nonumber
\\ & &
\Rightarrow 2(p-1)! + 1  \equiv  0 \pmod{p+2} \qquad \textrm{ iff
$p+2$ is odd prime.} \nonumber \eseq \bseq & & \Rightarrow
2[2(p-1)! + 1]  \equiv  0 \pmod{p+2} \qquad \textrm{ iff $p+2$ is
odd prime.} \nonumber
\\ & &
\Rightarrow 2[2(p-1)!+1]+(p+2)  \equiv  0 \pmod{p+2} \qquad
\textrm{ iff $p+2$ is odd prime.} \nonumber
\\ & &
\Rightarrow 4[(p-1)!+1]+p  \equiv  0 \pmod{p+2} \qquad \textrm{
iff $p+2$ is odd prime.} \nonumber
\\ & &
\Rightarrow 4[(p-1)!+1]  \equiv  -p \pmod{p+2} \qquad \textrm{ iff
$p+2$ is odd prime.} \eseq From (1) and (2), \bseq & & 4[(p-1)!+1]
\equiv -p \pmod{p(p+2)}   \qquad \textrm{ iff $p$, $p+2$ are odd
prime numbers.} \qquad \square \nonumber \eseq

\section{Derivation from Wilson's Theorem for Polignac
Conjecture} In the derivation of Clement's result, a method of
determining unknown coefficients has been applied. This method is
used to combine two congruence identities into one. Suppose we
have
\bseq f_1(p)+ C_1 & \equiv & 0 \pmod{p} \nonumber\\
\textrm{and } \qquad g(\lambda)f_1(p)+ C_2 & \equiv & 0 \pmod{p+
\lambda},\nonumber \eseq  where $(p, p+ \lambda)=1$ and
$(g(\lambda), p(p+\lambda))=1$, we establish two congruences:
\bseq X g(\lambda)[f_1(p)+C_1] + Y p & \equiv & 0 \pmod{p} \\
\textrm{and }\qquad X [g(\lambda)f_1(p)+C_2]+Y(p+ \lambda) &
\equiv & 0 \pmod{p+ \lambda}. \eseq If we choose X and Y in such a
way that $(X,p)=1$ and $(X, p+ \lambda)=1$, and \beq X
g(\lambda)C_1=X C_2+Y \lambda \nonumber \eeq the (3) and (4) is
equivalent to \beq X g(\lambda)f_1(p)+Y p+ X C_2 +Y\lambda \equiv
0 \pmod{p(p+ \lambda)} . \eeq

Specifically, in deriving Clement's result, we used $X=2$ and
$Y=1$. This method can be employed again in attempting the
Polignac Conjecture which states that for every positive number
$k$, there are infinitely many pairs of primes in the form of $p$
and $p+2k$. We propose a necessary and sufficient condition when
$p$ and $p+2k$ are primes.
\\ \\ \textbf{Theorem 3} \qquad Suppose $p \nmid 2k(2k)!$ and p is
odd. Then \bseq 2k(2k)![(p-1)!+1] \equiv [1-(2k)!]p \pmod{p(p+2k)}
& & \textrm{iff p, p+2k are odd primes.} \nonumber \eseq
\textbf{Proof.} \qquad From Wilson's theorem, we have \bseq
(p-1)! +1 & \equiv & 0 \pmod{p} \textrm{ iff $p$ is odd prime.} \nonumber \\
(p+2k-1)! +1 & \equiv & 0 \pmod{p+2k} \textrm{ iff $p+2k$ is odd
prime.}\nonumber \eseq The LHS of the congruence can be rewritten
as \bseq(p+2k-1)(P+2k-2)\cdots p (p-1)! +1 & \equiv & (-1)(-2)
\cdots (-2k)(p-1)!+1 \nonumber \\ {} & \equiv & (2k)!(p-1)!+1
\equiv  0 \pmod{p+2k} . \nonumber \eseq By setting $\lambda = 2k$,
$g(\lambda)=(2k)!$, $f_1(p)=(p-1)!$, $C_1=C_2=1$, $X=2k$ and
$Y=(2k)!-1$ in (5), we obtain \beq 2k(2k)![(p-1)!+1] \equiv
[1-(2k)!]p \pmod{p(p+2k)} \quad \nonumber \eeq \qquad \qquad
\qquad iff p, p+2k are odd primes, assuming $p \nmid 2k(2k)!$ and
p is odd. \qquad \qquad $\square$
\\
\\For example, when $k=2$, we have, \beq \textrm{for $p\nmid96$,}
\quad 96[(p-1)!+1]\equiv -23p \pmod{p(p+4)}\quad \textrm{iff $p$,
$p+4$ are odd primes.} \nonumber \eeq

\section{Simionov's Generalization of Wilson's Theorem}
By elementary mathematical manipulations, we can generalize
Wilson's Theorem. Simionov had found out the following
generalization of Wilson's Theorem \cite{Florentin Smarandache}
\bseq \textrm{\textbf{Theorem 4}}\qquad (k-1)!(p-k)! \equiv (-1)^k
\pmod{p}, \quad \textrm{iff $p$ is prime}\quad \forall k\in
\mathbb{Z^+} , \quad k\leq p \quad \eseq \textbf{Proof.} By
Wilson's Theorem, \bseq  & & (p-1)!  \equiv  -1 \pmod{p} \textrm{
iff $p$ is prime.} \nonumber \\ & &
 \Leftrightarrow
(p-1)(p-2)(p-3)\cdots (p-k+1)(p-k)!  \equiv  -1 \pmod{p} \textrm{
iff $p$ is prime.}\nonumber \\ & &
 \Leftrightarrow
(-1)(-2)(-3)\cdots (-k+1)(p-k)!  \equiv  -1 \pmod{p} \textrm{ iff
$p$ is prime.}\nonumber \eseq \bseq & & \Leftrightarrow
(-1)^{k-1}(1)(2)(3)\cdots (k-1)(p-k)!  \equiv -1 \pmod{p} \textrm{
iff $p$ is prime.}\nonumber
\\ & &
\Leftrightarrow (-1)^{k-1}(-1)^{k-1}(1)(2)(3)\cdots (k-1)(p-k)!
\equiv  (-1)^{k} \pmod{p} \textrm{ iff $p$ is prime.}\nonumber
\\ & &
\Leftrightarrow (-1)^{k-1}(-1)^{k-1}(k-1)!(p-k)!  \equiv (-1)^{k}
\pmod{p} \textrm{ iff $p$ is prime.}\nonumber
\\ & &
\Leftrightarrow (k-1)!(p-k)!  \equiv  (-1)^{k} \pmod{p} \textrm{
iff $p$ is prime.} \nonumber \qquad \square\eseq Besides
Simionov's result, other generalizations include \\ \qquad
\textbf{Corollary 1}\qquad For two distinct odd prime numbers,
$p_{1}, p_{2}$ and a positive odd integer $p$ such that $p + 1 =
p_{1} + p_{2}$ (Partial consequence of Goldbach's Conjecture \cite
{HALE UNSAL}) \beq (p-p_{1})! (p-p_{2})! \equiv -1 \pmod{p}
\textrm{ iff $p$ is prime}\nonumber \eeq \textbf{Corollary
2}\qquad For integers $k_{1}$ and $k_{2}$ such that $0\leq k_{1} <
k_{2} < p$ and $k_{2} - k_{1} \equiv 1 \pmod{2}$, \beq k_{1}!
k_{2}! (p-k_{1}-1)! (p-k_{2}-1)! \equiv -1 \pmod{p} \textrm{ if
$p$ is prime} \nonumber \eeq However,these generalizations are
less significant compared to Sinionov's
Theorem. \\
\\
A direct consequence of Simionov's Theorem can be obtained by
substituting $k=\frac{p+1}{2}$ into the equation (6), \bseq
\textrm{\textbf{Theorem 5}} \qquad \qquad \qquad \qquad \quad
[(\frac{p-1}{2})!]^{2} & \equiv & (-1)^{\frac{p+1}{2}} \textrm{
iff $p$ is odd prime}\qquad \qquad \qquad \qquad \eseq
\textbf{Proof.} \quad This result can also be directly derived
from Wilson's Theorem as shown below, for an odd integer $p$,
\bseq & & (p-1)(p-2)(p-3) \cdots (\frac{p+1}{2})(\frac{p-1}{2})
\cdots 3\cdot 2\cdot 1  \equiv  -1 \pmod{p} \textrm{ iff $p$ is
odd prime.} \nonumber
\\ & &
\Leftrightarrow [(p-1)\cdot 1][(p-2)\cdot 2][(p-3)\cdot 3] \cdots
[(\frac{p+1}{2})(\frac{p-1}{2})]  \equiv  -1 \pmod{p} \textrm{ iff
$p$ is odd prime.} \nonumber
\\ & &
\Leftrightarrow (-1)(-4)(-9) \cdots (-(\frac{p-1}{2})^{2}) \equiv
-1 \pmod{p} \textrm{ iff $p$ is odd prime.} \nonumber
\\
& & \Leftrightarrow (-1^{2})(-2^{2})(-3^{2}) \cdots
[-(\frac{p-1}{2})^{2}]  \equiv  -1 \pmod{p} \textrm{ iff $p$ is
odd prime.} \nonumber
\\ & & \Leftrightarrow
(-1)^{\frac{p-1}{2}}[(\frac{p-1}{2})!]^{2}  \equiv  -1 \pmod{p}
\textrm{ iff $p$ is odd prime.} \nonumber
\\ & & \Leftrightarrow
(-1)^{\frac{p-1}{2}}(-1)^{\frac{p-1}{2}}[(\frac{p-1}{2})!]^{2}
\equiv (-1)(-1)^{\frac{p-1}{2}} \pmod{p} \textrm{ iff $p$ is odd
prime.} \qquad \nonumber
\\ & &
\Leftrightarrow [(\frac{p-1}{2})!]^{2}   \equiv
(-1)^{\frac{p+1}{2}} \pmod{p} \textrm{ iff $p$ is odd prime.}
\qquad \nonumber \square \eseq As the calculation for the
factorial of a number can be tremendous, this result can reduce
the computation for determining the primality of a number, as
$[(\frac{p-1}{2})!]^2$ is much smaller than $(p-1)!$, especially
for large $p$. Hence, Simionov's finding can be viewed as an
improvement of Wilson's Theorem. We can take this approach when
solving Twin Prime Conjecture and Polignac Conjecture.

\section{Derivation from Simionov's Theorem for Polignac Conjecture}

From Simionov's result, we use the method of the mathematical
manipulations mentioned in \textbf{Section III} to derive a
sufficient and necessary condition for twin prime numbers to
exist, i.e. for $p$ and $p+2$ to be primes \bseq
\textrm{\textbf{Theorem 6}} \qquad \qquad \quad & &
2[(\frac{p-1}{2})!]^2  + (-1)^{\frac{p-1}{2}} (5p+2) \equiv 0
\pmod{p(p+2)} \qquad \qquad \qquad \nonumber \\ & & \textrm{ iff
$p$ and $p+2$ are odd prime numbers.}  \qquad \qquad \nonumber
\eseq \textbf{Proof.} \qquad From (7), for odd integer $p$,\bseq
 & &
[(\frac{p-1}{2})!]^2+(-1)^{\frac{p-1}{2}}  \equiv  0 \pmod{p}
\textrm{ iff $p$ is odd prime.} \nonumber
\\ & &
\Rightarrow 2[(\frac{p-1}{2})!]^2+2(-1)^{\frac{p-1}{2}} \equiv  0
\pmod{p} \textrm{ iff $p$ is odd prime.} \nonumber
\\ & &
\Rightarrow 2[(\frac{p-1}{2})!]^{2}+(-1)^{\frac{p-1}{2}}(2+5p)
  \equiv  0 \pmod{p} \textrm{ iff $p$ is odd prime.}
 \\ & & \textrm{Similarly for
  $p+2$,}
  \nonumber \\  & &
[(\frac{p+1}{2})!]^2+(-1)^{\frac{p+1}{2}}  \equiv  0 \pmod{p+2}
\textrm{ iff $p+2$ is odd prime} \nonumber
\\ & &
\Rightarrow 8[(\frac{p+1}{2})!]^{2}+8(-1)^{\frac{p+1}{2}}
 \equiv  0 \pmod{p+2} \textrm{ iff $p+2$ is odd prime.}
\nonumber
\\ & &
\Rightarrow
8({\frac{p+1}{2}})^2[(\frac{p-1}{2})!]^{2}+(-8)(-1)^{\frac{p-1}{2}}
\equiv
 0 \pmod{p+2} \textrm{ iff $p+2$ is odd prime.} \nonumber
\\ & &
\Rightarrow
2(p^2+2p+1)[(\frac{p-1}{2})!]^{2}+[5(p+2)-8]^{\frac{p-1}{2}}
\equiv 0 \pmod{p+2} \textrm{ iff $p+2$ is odd prime.} \nonumber
\\ & &
\Rightarrow 2[(\frac{p-1}{2})!]^{2}+(5p+2)(-1)^{\frac{p-1}{2}}
 \equiv  0 \pmod{p+2} \textrm{ iff $p+2$ is odd prime.}
 \eseq
\qquad \qquad From (8) and (9), \beq
2[(\frac{p-1}{2})!]^2+(5p+2)(-1)^{\frac{p-1}{2}} \equiv  0
\pmod{p(p+2)} \textrm{ iff $p$, $p+2$ are odd primes.} \nonumber
\quad \square \eeq  If we use the same method again on the pair of
congruences \bseq & & [(\frac{p-1}{2})!]^2+(-1)^{\frac{p-1}{2}}
\equiv  0 \pmod{p} \textrm{ iff $p$ is odd prime} \nonumber \\ & &
[(\frac{p+2k-1}{2})!]^2+(-1)^{\frac{p+2k-1}{2}}  \equiv  0
\pmod{p+2k} \textrm{ iff $p+2k$ is odd prime} \nonumber \eseq
 and setting $X=4^k (2k)$ and $Y=(2k-1)!!^2 (-1)^{\frac{p-1}{2}}- 4^k
 (-1)^{\frac{p+2k-1}{2}}$ in (5) where $(2n-1)!!$ denotes $\prod_{i=1}^n (2i-1)$, it easily
 follows that \\ \textbf{Theorem 7} \qquad Suppose $p\nmid (2k-1)!!^2$ and $p$ is
odd, \bseq 2k(2k-1)!!^2[(\frac{p-1}{2})!]^2+
(-1)^{\frac{p-1}{2}}[(2k-1)!!^2(p+2k)+4^k(-1)^{k+1}p] &  \equiv 0
& \pmod{p(p+2k)} \quad \nonumber \\ \textrm{iff $p$, $p+2k$ are
primes.} & & \nonumber \eseq This is a necessary and sufficient
condition for $p$ and $p+2k$ to be primes assuming $p\nmid
(2k-1)!!^2$ and $p$ is odd. For example, when $k=2$, we have,
supposing $p\nmid 9$ and $p$ is odd, \beq 36[(\frac{p-1}{2})!]^2+
(-1)^{\frac{p-1}{2}}[36-7p] \equiv 0 \pmod{p(p+4)}\quad
\textrm{iff $p$, $p+4$ are primes.}\nonumber \eeq

\section{Conclusion}
In this paper, we have proposed necessary and sufficient
conditions for $p$ and $p+2k$ to be primes where $k\in
\mathbb{\mathbb{Z}^{+}}$(Theorem 3 and Theorem 7). Unfortunately,
because of the tedious calculation of the factorial of a large
number, these theorems may not be very effective in searching for
large twin primes. Though Theorem 7 reduces significantly the
number of computational steps of Theorem 3, it is still a daunting
task to calculate $(\frac{p-1}{2})!$ when p becomes very large.
Moreover, as we cannot prove that there are infinitely many
numbers satisfying the condition yet, we cannot solve the Polignac
Conjecture. Thus, we need to further improve these theorems or
look for other alternative routes in attempting the
Twin Prime Conjecture and the Polignac Conjecture. \\ \\
\\\textbf{ ACKNOWLEDGEMENT } \\ We would like to express our
gratitude to Assoc Prof Choi Kwok Pui of National University of
Singapore in guiding us in the project on Benford's Law which led
us to examine the goodness-of-fit for prime number distribution to
the Benford's distribution. While tackling this intriguing
problem, we extended our knowledge of primes and obtained the
results stated in this paper along the way.

\end{document}